\documentclass[a4paper,16pt]{article}
\usepackage{amsthm}
\usepackage[sumlimits]{amsmath}
\usepackage{amsfonts}
\usepackage{mathrsfs}
\usepackage{varioref}
\usepackage{apalike}
\usepackage{caption2}
\usepackage{colordvi,graphics,latexsym}
\allowdisplaybreaks

\def\MM#1{\boldsymbol{#1}}
\newtheorem{theorem}{Theorem}[section]

\DeclareMathOperator{\diff}{d}

\newcommand{\pp}[2]{\frac{\partial #1}{\partial #2}} 
\newcommand{\dd}[2]{\frac{\delta #1}{\delta #2}}

\newtheorem{definition}[theorem]{Definition}

\begin{document}
\sffamily

\title{A General Approach for Producing Hamiltonian Numerical Schemes for Fluid
  Equations}
\author{Colin Cotter}
\maketitle

\begin{abstract}
  Given a fluid equation with reduced Lagrangian $l$ which is a
  functional of velocity $\MM{u}$ and advected density $D$ given in
  Eulerian coordinates, we give a general method for semidiscretising
  the equations to give a canonical Hamiltonian system; this system
  may then be integrated in time using a symplectic integrator. The
  method is Lagrangian, with the variables being a set of Lagrangian
  particle positions and their associated momenta. The canonical
  equations obtained yield a discrete form of Euler-Poincar\'e
  equations for $l$ when projected onto the grid, with a new form of
  discrete calculus to represent the gradient and divergence
  operators. Practical symplectic time integrators are suggested for a
  large family of equations which include the shallow-water equations,
  the EP-Diff equations and the 3D compressible Euler equations, and
  we illustrate the technique by showing results from a numerical
  experiment for the EP-Diff equations.
\end{abstract}

\tableofcontents

\section{Introduction}
There has been much recent interest in obtaining Hamiltonian methods
for the various equations of motion for fluids (such as the
shallow-water equations, and the 2D Euler equations). The Hamiltonian
method programme consists of two stages:
\begin{enumerate}
\item
  Take a Hamiltonian PDE with Hamiltonian $\mathcal{H}$ and structure
  operator $\mathcal{J}$ so that the equation takes the form
  $$
  \mathcal{J}\MM{z}_t(\MM{x}) = \dd{\mathcal{H}}{\MM{z}}.
  $$
  Discretise the Hamiltonian and the structure operator to obtain a
  finite dimensional Hamiltonian $\hat{H}$ and symplectic structure
  matrix $\hat{J}$ to give a finite dimensional Hamiltonian system
  $$
  \hat{J}\dot{\MM{z}} = \pp{\hat{H}}{\MM{z}}.
  $$
\item Discretise the finite dimensional system in time using a
  symplectic integrator (for a general review of symplectic methods
  see \cite{hair02:geometric_numerical_integration} or
  \cite{leim05:sim_ham_dyn}).
\end{enumerate}
The advantage of this approach is that the symplectic integrator will
guarantee excellent long-time conservation properties
\cite{ben94:_ham_int_near_iden_sym_map_app_sym,hair94:bac,HaiLub97,reic99:bac}
with the spatially-discrete Hamiltonian $\hat{H}$ satisfying
$$
|\hat{H}(t)-\hat{H}(0)|<c_1\delta t^p, \qquad |t|<c_2e^{-c_3/\delta t},
$$
where $p$ is the order of the symplectic method in time and where
$c_1$, $c_2$ and $c_3$ are positive constants. Note that Hamiltonian
methods are not the same as multi-symplectic integrators
\cite{j.98:mul_pdes,mar01,reic00:fin_pdes,moor03:mul_int_met_ham_pdes,moor03:bac}
where the Lagrangian is symmetrically discretised in time and space,
and discrete variations are taken to obtain the numerical method,
which then offers additional symmetry properties
\cite{oliv04:app}.

In general it has been very difficult to obtain discrete structure
operators $\hat{J}$ which satisfy the Jacobi identity and progress has
only been made for Hamiltonian PDEs with canonical structure (with a
few exceptions \emph{e.g.} \cite{zeit91:_fin_anal_ideal_hyd}). In the
context of fluid equations this means that the equations must be
solved in Lagrangian coordinates where the symplectic structure is
canonical. For the case of 2D incompressible flow this leads to
methods where Lagrangian particles carry vorticity instead of mass,
either as Dirac $\delta$-functions (point-vortex methods \cite[for a
review]{cot00:vortex_methods}) or in ``blobs'' surrounding the
particle (vortex-blob methods
\cite{chor73:num,oliv01:new,cot04:geom}).  For the case of
compressible flow the particles carry masses which are interpolated
through basis functions to the entire domain to give the density
field; this type of discretisation is called Smoothed-Particle
Hydrodynamics (SPH) \cite{gingold77,lucy77,sal83:prac}. The related
Hamiltonian Particle-Mesh (HPM) method
\cite{fran02:hamiltonian,cot03:hamiltonian_particle_mesh_method_two,fran04:ham_par_mes,cot04:md_nw}
discretises a regularised form of the equations (with the
regularisation operator discretised on an Eulerian grid) so that good
long-time behaviour is obtained when the finite dimensional system is
integrated with a symplectic integrator (this is described further in
section \ref{background}).

The difficulty with Lagrangian coordinates is that the equations can
become very complicated (especially when differential operators in
Eulerian coordinates are involved). Analytically, the most tractable
approach is to transform the Lagrangian $L$ for the flow maps in
Lagrangian coordinates to a reduced Lagrangian $l$ for the velocity
field $\MM{u}(\MM{x})$ in Eulerian coordinates, and then to use the
Euler-Poincar\'e (EP) equations
\cite{hol98:eul_poin,hol98:eul_poin_equat_geop_fluid_dyn} for $\MM{u}$
which are equivalent to the canonical Hamilton's equations for the
Lagrangian variables.  Writing the reduced Lagrangian which integrates
over Eulerian coordinates means that the equations are much easier to
handle. This greatly aids model derivation
\cite{hol02:eul_poin_dyn_per_com_fluid}, analysis
\cite{mar00:eul,mar01:lag_nav_stok_lan,shk02:lag_eul_lae}, model
reduction \cite{hol98,oliv05:var} and averaging
\cite{hol99:_fluc_lag_eul}; it is clear that it would be attractive to
use the reduced Lagrangian to derive numerical methods too.  However,
in the case of general fluid equations this would mean losing the
canonical structure operator because the Legendre transform in
Eulerian coordinates results in a non-canonical Lie-Poisson equation.

In this paper we give a general approach, using Lagrangian variables
so that there is a canonical structure operator which is easy to
discretise, whilst evaluating the Hamiltonian on an Eulerian
grid. This approach fully develops a dual grid-particle approach
(as an extension of the mixed form that HPM uses with the potential
energy in Eulerian variables whilst the kinetic energy is in
Lagrangian variables) The particle-mesh formulation allows quantities
(such as momentum and density) which are given at the particle
locations, to be interpolated to the Eulerian grid. The Hamiltonian is
then approximated as a sum over the Eulerian grid points, and Theorem
\ref{my theorem} in this paper shows that the canonical equations
obtained from the discrete Hamiltonian produce a set of discrete EP
equations on the grid (although computationally the equations are
treated in Lagrangian form). It then remains to apply a numerical
time-stepping algorithm which preserves the canonical symplectic
structure.

The structure of the rest of the paper is as follows: section
\ref{background} provides some background on SPH and HPM to prepare
the way for the general approach to obtaining semidiscretisations.
This approach is set out, together with the statement and proof of
Theorem \ref{my theorem}, in section \ref{approach}, and suitable
symplectic time-integrators are suggested in section \ref{time int}.
The next two sections show how to apply the approach in two examples:
the shallow-water-$\alpha$ equations in section \ref{swalpha}, and the
2D EP-Diff equation in section \ref{epdiff}. Both of these examples
are in 2-dimensional geometry although the approach is completely
general and may be used to discretise $n$-dimensional PDEs. Section
\ref{numerics} gives some numerical results for the 2D EP-Diff
equation and section \ref{discussion} contains some brief final
discussion.
 
\section{Background}
\label{background}
HPM was originally introduced in \cite{fran02:hamiltonian} as a
semidiscretisation of the regularised layer-depth shallow-water
(RLDSW) equations:
$$
\MM{u}_t+\MM{u}\cdot\nabla\MM{u}=-g\nabla(1-\alpha^2\Delta)^{-1}D, \qquad
D_t+\nabla\cdot(D\MM{u})=0.
$$
The success of this method is based on its ``dual picture'' nature,
as the Lagrangian particle dynamics provides a canonical symplectic
structure (which is easy to treat numerically using symplectic
integrators) whilst the inverse modified Helmholtz operator takes its
simplest form on the Eulerian grid.

HPM can be viewed as a modification of SPH \cite{gingold77,lucy77}
which is a Lagrangian particle method for compressible flow. In the
SPH method for shallow-water, the layer-depth $D$ is represented by a
linear combination of radial basis functions which are centred on a
finite set $\{\MM{X}_{\beta}(t)\}_{\beta=1}^n$ of Lagrangian
particles:
$$
D(\MM{x},t) = \frac{1}{\Delta S}\sum_{\beta=1}^n 
D_{\beta}\phi(\MM{x}-\MM{X}_{\beta}(t)),
$$ 
where $\phi$ is a radial basis function with integral $\Delta S$,
and $\{D_{\beta}\}_{\beta=1}^n$ are constant weights. Other fields
$f(\MM{x},t)$ (such as temperature, velocity \emph{etc.}) may be
interpolated from the values $\{f_{\beta}(t)\}_{\beta=1}^n$ at the
particle locations to the entire domain \emph{via} the product
$$
D(\MM{x},t)f(\MM{x},t)=\frac{1}{\Delta S}\sum_{\beta=1}^n
D_{\beta}f_{\beta}(t)\phi(\MM{x}-\MM{X}_{\beta}(t)).
$$
This interpolation may then be used to establish the continuity equation
for $h$:
\begin{eqnarray*}
\partial_tD(\MM{x},t) & = & \frac{1}{\Delta S}\partial_t\sum_{\beta=1}^n D_{\beta}\phi(\MM{x}-\MM{X}_{\beta}(t)), \\
& = & \frac{1}{\Delta S}\sum_{\beta=1}^nD_{\beta}\dot{\MM{X}}_{\beta}(t)\nabla_{\MM{X}_{\beta}}\phi(\MM{x}-\MM{X}_{\beta}(t)), \\
& = & -\frac{1}{\Delta S}\sum_{\beta=1}^nD_{\beta}\dot{\MM{X}}_{\beta}(t)\nabla_{\MM{x}}\phi(\MM{x}-\MM{X}_{\beta}(t)), \\
& = & -\nabla_{\MM{x}}\cdot\frac{1}{\Delta S}\sum_{\beta=1}^n D_{\beta}\dot{\MM{X}}_{\beta}(t)\phi(\MM{x}-\MM{X}_{\beta}(t)), \\
& = & -\nabla_{\MM{x}}\cdot(D\MM{v}(\MM{x},t)),
\end{eqnarray*}
where $\nabla_{\MM{X_{\beta}}}$ is the gradient with respect to $\MM{X}_{\beta}$,
$\nabla_{\MM{x}}$ is the gradient with respect to $\MM{x}$, and
$\MM{v}(\MM{x},t)$ is the interpolation of the particle velocities
$\dot{\MM{X}}_{\beta}$ to the entire domain.

In SPH methods, the forces on the particles are simply evaluated by
interpolating the forces to the entire domain and then substituting
the particle positions in the formula. These methods are thus
mesh-free which makes them popular with astrophysicists who are able
to model self-gravitating fluids (such as colliding stars and
accretion discs) contained in an infinite vacuum simply by giving the
basis functions finite support so that there is no density in the
far-field. In particle-mesh methods, however, the aim is to meet the
need for easy evaluation of Eulerian differential operators such as
the inverse modified Helmholtz operator in the RLDSW equations. For a
general mesh with $m$ grid points $\{\MM{x}_k\}_{k=1}^m$, the
layer-depth $D_k$ at grid point $k$ is given by
$$
D(\MM{x}_k,t) = \frac{1}{\Delta S}\sum_{\beta=1}^n D_\beta\psi_{k}(\MM{X}_{\beta}),
$$
where 
$$
\psi_{k}(\MM{X})=\phi(\MM{x}_k;\MM{X}),
$$
and $\phi(\MM{x};\MM{X}_{\beta})$ is some general basis function
(such as a radial basis function or a Cartesian product of 1D basis
functions) centred on $\MM{X}_{\beta}$ that satisfies
$$
\nabla_{\MM{x}}\phi(\MM{x};\MM{X})=-\nabla_{\MM{X}}\phi(\MM{x};\MM{X}),
\qquad \int \phi(\MM{x};\MM{0})\diff^2x=\Delta S.
$$ 
The key to the method is the
restriction to basis functions (such as the Cartesian product of cubic
B-splines \cite[for example]{boor78}) which satisfy the partition of
unity property:
$$
\sum_{\beta=1}^n\psi_{k}(\MM{x})=1 \qquad \forall \MM{x}
$$
which means that any quantity $f$ may be interpolated 
from the grid to the whole domain:
$$
f(\MM{x})=\sum_{k=1}^mf(\MM{x}_k)\psi_k(\MM{x}).
$$
This allows us to interpolate from the grid to the particles
\emph{via}
$$
f_{\beta} = f(\MM{X_{\beta}}) = \sum_{k=1}^mf(\MM{x}_k)\psi_k(\MM{X}_{\beta}).
$$

The gradient of the density can be obtained by taking the gradient of
the basis functions evaluated at the particle positions:
\begin{eqnarray*}
\nabla_{\MM{x}}
D(\MM{x})|_{\MM{x}=\MM{X}_{\beta}} &=& - \sum_{k=1}^mD(\MM{x}_k)\nabla_{\MM{X}_{\beta}}\psi_k(\MM{X}_{\beta}) \\
&=& - \sum_{k=1}^m\nabla_{\MM{X}_{\beta}}\psi_k(\MM{X}_{\beta})
\frac{1}{\Delta S}\sum_{\gamma=1}^n
D_\gamma\psi_k(\MM{X}_\gamma). \\
\end{eqnarray*}

The RLDSW equations were originally introduced after failed attempts
to exploit the canonical Hamiltonian structure of SPH (applied to the
SW equations) by performing the time integration using a symplectic
method (with the aim of obtaining long time approximate preservation
of energy, momentum and adiabatic invariants). These attempts failed
because computational efficiency necessitates the use of basis
functions with compact support. When the flow is nearly incompressible
(\emph{e.g.} when the system is near to geostrophic balance in the
rotating case) the basis functions supported local compressible
oscillations which quickly destabilised the flow and led to
equipartition of energy much too quickly. The smoothing operator
applied to the layer depth means that the local basis functions become
global and the small oscillations are filtered out in the momentum
equation. The semi-discretised equations are:
\begin{eqnarray*}
\dot{\MM{V}}_{\beta} & = & g\sum_{k=1}^m\nabla_{\MM{X}_{\beta}}\psi_k(\MM{X}_{\beta})
\sum_{l=1}^m(A^{-1})_{kl}\frac{1}{\Delta S}\sum_{\gamma=1}^n\psi_k(\MM{X}_\gamma)D_{\gamma}, \\
\dot{\MM{X}}_{\beta} & = & \MM{V}_\beta, \qquad \beta=1,\ldots,n,
\end{eqnarray*}
where $A$ is the numerical approximation to the modified Helmholtz
operator (obtained using \emph{e.g.} Fourier series or finite
elements). These equations have a canonical Hamiltonian structure
with Hamiltonian
$$
H =
\sum_{\beta=1}^n\frac{|\MM{m}_{\beta}|^2}{2D_{\beta}}+g\sum_{k,l=1}^m
\frac{1}{\Delta S}\sum_{\beta,\gamma=1}^n
\psi_k(\MM{X}_{\beta})(A^{-1})_{kl}D_{\gamma}\psi_l(\MM{X}_{\gamma}),
$$
where the momentum is given by
$$
\MM{m}_{\beta}= D_{\beta}\MM{v}_{\beta}, \qquad \beta=1,\ldots,n.
$$
These equations may be integrated using the symplectic
St\"ormer-Verlet method because the Hamiltonian splits into separate
functions of $\{\MM{m}\}_{\beta=1}^n$ and $\{\MM{X}\}_{\beta=1}^n$ (as
in classical mechanics) \cite{mcl99:odes}.

In this method, the potential energy is obtained by interpolating the
height field to the grid and summing over the gridpoints, whilst the
kinetic energy is evaluated in Lagrangian coordinates only. In the
next section we shall extend the approach to the kinetic energy used
here by evaluating the entire Hamiltonian on the Eulerian grid.

\section{Canonical Hamiltonian particle-mesh semidiscretisations}
\label{approach}
Throughout this section, and for the rest of this paper, we shall
adopt the notation that a bar indicates a quantity stored at a
particle location \emph{e.g.}  $\overline{\MM{m}}$, $\overline{D}$,
\emph{etc.}, and a tilde indicates a quantity stored at a grid
location \emph{e.g.}  $\tilde{\MM{m}}$, $\tilde{D}$, $\tilde{\MM{u}}$.
Also, Greek letters will be used for indices denoting the particle
numbering, \emph{e.g.} momentum $\overline{\MM{m}}_\beta$ at particle
$\MM{X}_\beta$, whilst Roman letters will be used for indices denoting
the grid numbering, \emph{e.g.} velocity $\tilde{\MM{u}}_k$ at grid
point $\MM{x}_k$.

The programme for producing Hamiltonian particle-mesh
semidiscretisations is as follows:
\begin{enumerate}
\item Start with the reduced Lagrangian $l$ for the EP
  equations given in Eulerian coordinates. This method applies to
  semi-direct product systems where the $l$ is a functional of the
  velocity $\MM{u}$ and the density $D$ (which satisfies the
  continuity equation).  In this section we shall look at PDEs in two
  dimensions, but the extension to higher (and lower) dimensions is
  straightforward.

\item To obtain the discrete Lagrangian $\hat{L}$, take the continuous
  Lagrangian $l$ which is written as an integral over $\mathcal{D}$
  $$
  l = \int \lambda(\MM{u},D)\diff{^2x},
  $$
  where $\lambda$ is a nonlinear operator, and replace it with a
  sum over gridpoints
  $$
  \hat{L} = \sum_{kl}M_{kl}\lambda_k\left(\{\tilde{\MM{u}}_l\}_{l=1}^m,
    \{\tilde{D}_l\}_{l=1}^m\right),
  $$
  where $M$ is the matrix representing the approximation to
  integration (often referred to in the finite-element literature as
  the mass matrix) and where $\lambda_k$ is a function approximating
  $\lambda$ at the point $\MM{x}_k$, with any differential operators
  being replaced by matrix operations on the grid in the manner
  consistent with the choice of $M$.

\item Use the Legendre transform 
  $$
  \hat{H}(\tilde{\MM{m}}_l\}_{l=1}^m,\{\tilde{D}_l\}_{l=1}^m))
  = \sum_{kl}M_{kl}\tilde{\MM{u}}_k\cdot\tilde{\MM{m}}_k
  -\hat{L}, 
  $$
  where
  \begin{equation}
    \label{discrete H}
  \sum_{l}M_{kl}\tilde{\MM{m}}_l=\pp{\hat{L}}{\tilde{\MM{u}}_k},
  \end{equation}
  to obtain the Hamiltonian written in Eulerian coordinates.

\item Take $n$ Lagrangian particles $\{\MM{X}_\beta\}_{\beta=1}^n$
  together with $m$ Eulerian grid points $\{\MM{x}_k\}_{k=1}^m$ on the
  domain $\mathcal{D}$ and a set of basis functions $\psi_k(\MM{x})$
  defined by
  $$
  \psi_k(\MM{x})=\phi(\MM{x}_k;\MM{x}),
  $$
  with the partition-of-unity property
  $$
  \sum_k\psi_k(\MM{x})=1\mbox{ for all }\MM{x}\in\mathcal{D}.
  $$
  and where $\phi$ satisfies
  $$
  \nabla_{\MM{x}}\phi(\MM{y};\MM{x})=-\nabla_{\MM{y}}\phi(\MM{y};\MM{x}), 
  \qquad
  \int_{\mathcal{D}}\phi(\MM{x};\MM{y})\diff^2x=\Delta S,
  $$
  
  This property allows grid values $f(\MM{x}_k)$ to be interpolated to
  the whole domain $\mathcal{D}$:
  $$
  f(\MM{x}) = \sum_kf(\MM{x}_k)\psi_k(\MM{x}),
  $$
  as well as gradients:
  $$
  \nabla f(\MM{x}) = \sum_kf(\MM{x}_k)\nabla\psi_k(\MM{x}).
  $$
  Here we avoid the issue of boundary conditions by stipulating
  that $\mathcal{D}$ is a closed, compact manifold (\emph{e.g.} a
  square with periodic boundaries).
  
\item Give the $\beta$-th particle a momentum density
  $\overline{\MM{m}}_\beta$ and a constant layer-depth density
  $\overline{D}_\beta$ (which will be identified as a conserved mass
  in the particle system that is obtained after semidiscretisation).
  Interpolate the layer-depth density and momentum density to the grid
  using
  \begin{equation}
    \label{interp dens}
    \sum_{l}M_{kl}\tilde{\MM{m}}_l =
    \sum_\beta\frac{\overline{\MM{m}}_\beta}{\Delta S}
    \psi_k(\MM{X}_\beta),\qquad
    \sum_{l}M_{kl}\tilde{D}_l = \sum_\beta\frac{\overline{D}_\beta}{\Delta S}
    \psi_k(\MM{X}_\beta).
  \end{equation}
  The mass matrix is present because the product of the momentum
  $\overline{\MM{m}}$ with the basis functions $\psi$ should be
  interpreted as the approximation to the one-form density
  $\MM{m}\diff^2x$; inverting the mass matrix ``strips off'' the
  density $\diff^2x$.

\item The equations of motion for the particle positions
  $\{\MM{X}_\beta\}_{\beta=1}^n$ and their associated local momenta
  $\{\overline{\MM{m}}_\beta/\Delta S\}_{\beta=1}^n$ are then the
  canonical Hamilton's equations for this discrete Hamiltonian
  (\ref{discrete H}) having applied the substitution in equations
  (\ref{interp dens}).
\end{enumerate}

Now we shall show that the equations of motion represent an
approximation to the EP equations. First, we make a few definitions
which set out our discrete calculus which is used to form the discrete
EP equations.
\begin{definition}
  For a quantity $\tilde{f}$ specified on the grid, define the
  \emph{grid-to-particle map} by interpolating to the whole of
  $\mathcal{D}$,
  $$
  f(\MM{x}) = \sum_k\tilde{f}_k\psi_k(\MM{x}),
  $$ 
  and evaluating at the particle positions
  $$
  f(\MM{X}_i) = \sum_k\tilde{f}_k\psi_k(\MM{X_i}).
  $$
  We shall use the following notation
  $$
  \left[\tilde{f}\right]_\beta = \sum_k\tilde{f}_k\psi(\MM{X}_\beta).
  $$
  Further define the \emph{grid-to-particle gradient map} applied
  to $\tilde{f}$ by interpolating to the whole of $\mathcal{D}$,
  $$
  f(\MM{x}) = \sum_k\tilde{f}_k\psi_k(\MM{x}),
  $$
  taking a gradient,
  $$
  \nabla f(\MM{x}) = \sum_k\tilde{f}_k\nabla\psi_k(\MM{x}),
  $$
  and evaluating at the particle positions
  $$
  \nabla f(\MM{X}_\beta) =
  \sum_k\tilde{f}_k\nabla\psi_k(\MM{X_\beta}).
  $$
  We shall use the following notation
  $$
  \left[\nabla\tilde{f}\right]_\beta =
  \sum_k\tilde{f}_k\nabla\psi_k(\MM{X}_\beta).
  $$
\end{definition}
\begin{definition}
  For a density $\overline{f}$ specified on the particle
  locations, define the \emph{particle-to-grid map} that interpolates
  to the grid by
  $$
  \sum_{l}M_{kl}\langle\overline{f}\rangle_l=
  \sum_\beta\frac{\overline{f}_\beta}{\Delta S}\psi_k(\MM{X}_\beta).
  $$ 
  Also for a function $\tilde{g}$ specified
  on the grid, use the following notation for the 
  interpolated product:
  $$
  \sum_{l}M_{kl}\langle\overline{f}[\tilde{g}]\rangle_l=
  \sum_\beta\frac{\overline{f}_\beta}{\Delta S}
  [\tilde{g}]_\beta\psi_k(\MM{X}_\beta).
  $$ 
  For a one-form density $\overline{\MM{f}}$ specified on the particle
  locations, further define the \emph{particle-to-grid divergence map} on
  the grid by
  $$
  \sum_{l}M_{kl}\langle\nabla\cdot\overline{\MM{f}}\rangle_l = 
  -\sum_\beta\frac{\overline{\MM{f}}_\beta}{\Delta S}
  \nabla\psi_k(\MM{X_\beta}).
  $$
  This map is the adjoint to the \emph{grid-to-particle gradient map}.
\end{definition}
These interpolation maps represent the whole key to the approach: they
allow the dynamical equations to be interpreted either on the grid or
on the particle positions. In particular, the canonical symplectic
structure is used to derive the particle dynamics (and so most of the
focus computationally will be on the particle variables) whilst the
Eulerian interpretation makes it possible to solve for the velocity
given the momentum on the grid. The following theorem shows that the 
equations on the grid can be interpreted as a discrete EP equation:
\begin{theorem}
\label{my theorem}
  On the grid the equations of motion obtained from this procedure
  take the following form
  \begin{eqnarray}
    \label{grid EP m}
    \frac{d}{dt}\langle\overline{\MM{m}}\rangle_k
    +\left\langle\nabla\cdot\left(\left[\tilde{\MM{u}}\right]
        \overline{\MM{m}}\right)\right\rangle_k
+\left\langle\left[(\nabla\tilde{\MM{u}})^T\right]
\cdot    
\overline{\MM{m}}\right\rangle_k
    &=&\left\langle \overline{D}
      \left[M^{-1}
        \nabla\pp{\hat{L}}{\tilde{D}}\right]\right\rangle_k, \\
    \label{grid EP D}
    \partial_t\langle\overline{D}\rangle_k+
    \left\langle \nabla\cdot \left(\overline{D}
        \left[\tilde{\MM{u}}\right]
      \right)\right\rangle_k
    & = & 0, \\
    \label{grid EP M defn}
    \langle\overline{\MM{m}}\rangle_k & = & 
\sum_l(M^{-1})_{kl}\pp{\hat{L}}{\tilde{\MM{u}}_l},
  \end{eqnarray}
  which gives an approximation on the grid, using the discrete
  calculus defined above, to the EP equations:
  \begin{eqnarray*}
  \partial_t\MM{m}+\nabla\cdot(\MM{u}\MM{m})+(\nabla\MM{u})^T\cdot\MM{m}
  &=& D\nabla\dd{l}{D}
  , \\
  \partial_tD+\nabla\cdot(\MM{u}D) &= &0,  \\
  \MM{m} & = & \dd{l}{\MM{u}}.
  \end{eqnarray*}
  Furthermore, at the particle locations the equations takes
  the form
  \begin{eqnarray}
    \label{particle EP m}
    \frac{d}{dt}\frac{\overline{\MM{m}}_\beta}{\overline{D}_\beta}
    + [(\nabla\tilde{\MM{u}})^T]_\beta\cdot
    \frac{\overline{\MM{m}}_\beta}{\overline{D}_\beta}
    &=&  
    \left[M^{-1}\nabla\pp{\hat{L}}{\tilde{D}}\right]_\beta, \\
    \label{particle EP X}
    \frac{d}{dt}\MM{X}_\beta &=& [\tilde{\MM{u}}]_\beta,
  \end{eqnarray}
  which gives an approximation at the particle locations to the EP
  equation in vector form:
  \begin{eqnarray*}
  \frac{D}{Dt}\frac{\MM{m}}{D}+(\nabla\MM{u})^T\cdot\frac{\MM{m}}{D}
  &=& \nabla\dd{L}{D}, \\
  \frac{D}{Dt} &=& (\partial_t+\MM{u}\cdot\nabla).
  \end{eqnarray*}
\end{theorem}
\emph{Proof.} The $\MM{X}$ equation obtained is 
  \begin{eqnarray*}
    \dot{\MM{X}}_\beta &=& \Delta S
    \nabla_{\overline{\MM{m}}_\beta}{\hat{H}}, \\
    & = & \sum_k\left(
      (\nabla_{\overline{\MM{m}}_\beta}\sum_lM_{kl}\tilde{\MM{m}}_l)\cdot\tilde{\MM{u}}_k +(\nabla_{\overline{\MM{m}}_\beta}\tilde{\MM{u}}_k)
      \cdot\sum_lM_{kl}\tilde{\MM{m}}_k\right) -
      \nabla_{\overline{\MM{m}}_\beta}\hat{L}.
  \end{eqnarray*}
  These three terms may be expanded out:
  \begin{eqnarray*}
    \sum_{k}\left(\nabla_{\overline{\MM{m}}_\beta}
      \sum_lM_{kl}\tilde{\MM{m}}_k\right)\cdot
    \tilde{\MM{u}}_k & = &
    \sum_k\left(\nabla_{\overline{\MM{m}}_\beta}
      \sum_\gamma\overline{\MM{m}}_\gamma
      \psi_k(\MM{X}_\gamma)\right)\cdot
    \tilde{\MM{u}}_k \\ 
    & = &
    \sum_k\tilde{\MM{u}}_k\psi_k(\MM{X}_\beta), \\
    \sum_k(\nabla_{\overline{\MM{m}}_\beta}\tilde{\MM{u}}_k)\cdot
    \sum_lM_{kl}\tilde{\MM{m}}_l & = &
    \sum_{klm}(\nabla_{\tilde{\MM{m}}_l}\tilde{\MM{u}}_k)\cdot
    M_{km}\tilde{\MM{m}}_m\psi_l(\MM{X}_\beta), \\
    \nabla_{\overline{\MM{m}}_\beta}\hat{L}
    & = & 
    \sum_{kl}    
(\nabla_{\tilde{\MM{m}}_l}\tilde{\MM{u}}_k)\cdot
(\nabla_{\tilde{\MM{u}}_k}\hat{L})\psi_l(\MM{X}_\beta)
    , \\
    & = & \sum_{klm}(\nabla_{\tilde{\MM{m}}_l}\tilde{\MM{u}}_k)
    \cdot M_{kl}\tilde{\MM{m}}_m
    \psi_l(\MM{X}_\beta),
\end{eqnarray*}
where we have made use of 
$$
\nabla_{\tilde{\MM{u}}_k}\hat{L}=\sum_lM_{kl}\tilde{\MM{m}}_l.
$$

This means that the second and third terms cancel and the $\MM{X}$
equation becomes
$$
\dot{\MM{X}}_\beta = \sum_k\tilde{\MM{u}}_k\psi_k(\MM{X}_\beta)=
[\tilde{\MM{u}}]_\beta.
$$

The $\overline{\MM{m}}$ equation is
\begin{eqnarray*}
\frac{\dot{\overline{\MM{m}}}_\beta}{\Delta S} &=& -
\nabla_{\MM{X}_\beta}\hat{H}
, \\
&=& \sum_{kl}\left(
(\nabla_{\MM{X}_\beta}M_{kl}\tilde{\MM{m}}_l)
\cdot\tilde{\MM{u}}_k+
(\nabla_{\MM{X}_\beta}\tilde{\MM{u}}_k)
\cdot M_{kl}\tilde{\MM{m}}_l\right)
-\nabla_{\MM{X}_\beta}\hat{L}.
\end{eqnarray*}
Once again, expand these three terms
\begin{eqnarray*}
\sum_{kl}(\nabla_{\MM{X}_\beta}M_{kl}\tilde{\MM{m}}_l)
\cdot\tilde{\MM{u}}_k & = &
\sum_k\tilde{\MM{u}}_k\cdot
\frac{\overline{\MM{m}}_\beta}{\Delta S}
\nabla_{\MM{X}_\beta}\psi_k(\MM{X}_\beta), \\
\sum_{kl}(\nabla_{\MM{X}_\beta}\tilde{\MM{u}}_k)
\cdot M_{kl}\tilde{\MM{m}}_k & = &
\sum_{klm}
\left(\frac{\overline{\MM{m}}_\beta}{\Delta S}\cdot
(\nabla_{\tilde{\MM{m}}_l}\tilde{\MM{u}}_k)\cdot
M_{km}\tilde{\MM{m}}_m\right)
\nabla_{\MM{X}_\beta}\psi_l(\MM{X}_\beta), \\
\nabla_{\MM{X}_\beta}\hat{L}
& = & 
\sum_{kl}
\left(\frac{\overline{\MM{m}}_\beta}{\Delta S}\cdot
(\nabla_{\tilde{\MM{m}}_l}\MM{u}_k)\cdot
\nabla_{\tilde{\MM{u}}_k}\hat{L}\right)
\nabla_{\MM{X}_\beta}\psi_l(\MM{X}_\beta) \\
& & \qquad -\sum_{k}\pp{\hat{L}}{\tilde{D}_k}
\nabla_{\MM{X}_\beta}\tilde{D}_k, \\
& = &
\sum_{klm}
\left(\frac{\overline{\MM{m}}_\beta}{\Delta S}\cdot
(\nabla_{\tilde{\MM{m}}_l}\tilde{\MM{u}}_k)\cdot
M_{km}\tilde{\MM{m}}_m\right)
\nabla_{\MM{X}_\beta}\psi_l(\MM{X}_\beta) \\
& & \qquad -\frac{\overline{D}_\beta}{\Delta S}
\sum_{kl}\pp{\hat{L}}{\tilde{D}_k}(M^{-1})_{kl}
\nabla_{\MM{X}_\beta}\psi_l(\MM{X}_\beta).
\end{eqnarray*}
The $\MM{m}$ equation may then be written in the form
$$
\frac{d}{dt}\overline{\MM{m}}_\beta = 
- 
[(\nabla\tilde{\MM{u}})^T]_\beta\cdot\overline{\MM{m}}_\beta
+\overline{D}_\beta\left[\nabla M^{-1}\pp{\hat{L}}{\tilde{D}}\right]_\beta,
$$ 
and dividing through by $\overline{D}_\beta$ gives equations
(\ref{particle EP m}-\ref{particle EP X}).

To obtain equation (\ref{grid EP m}), take the
time derivative of $\langle\overline{\MM{m}}\rangle_k$:
\begin{eqnarray*}
\frac{d}{dt}\sum_lM_{kl}\langle\overline{\MM{m}}\rangle_l & = & 
  \sum_k\left(\frac{d}{dt}\frac{\overline{\MM{m}}_\beta}{\Delta S}
  \right)\psi_k(\MM{X}_\beta)
+ \frac{1}{(\Delta x)^2}
\sum_\beta\frac{\overline{\MM{m}}_\beta}{\Delta S}
\frac{d}{dt}\psi_k(\MM{X}_\beta), \\
& = & \sum_lM_{kl}\left\langle\frac{d}{dt}\overline{\MM{m}}\right\rangle_l
+ 
\sum_\beta\frac{\overline{\MM{m}}_\beta}{\Delta S}
\frac{d}{dt}{\MM{X}}_\beta\cdot\nabla\psi_k
(\MM{X}_\beta), \\
& = & -\sum_lM_{kl}\Bigg(\left\langle
[(\nabla\tilde{\MM{u}})^T]\cdot\overline{\MM{m}}
\right\rangle_l+
\left\langle\overline{D}\left[
\nabla M^{-1}\pp{\hat{L}}{\tilde{D}}\right]\right\rangle_l
- \\
& & \qquad
\left\langle\nabla\cdot(\overline{\MM{m}}[\tilde{\MM{u}}]
)\right\rangle_l\Bigg),
\end{eqnarray*}
and dividing by $M$ gives the result. Similarly, equation (\ref{grid
  EP D}) follows by taking the time derivative of
$\langle\overline{D}\rangle_k$:
\begin{eqnarray*}
\frac{d}{dt}\sum_lM_{kl}\langle\overline{D}\rangle_l & = & 
\sum_\beta
\frac{\overline{D}_\beta}{\Delta S}
\frac{d}{dt}\psi_k(\MM{X}_\beta), \\
& = & 
\sum_\beta
\frac{\overline{D}_\beta}{\Delta S}\dot{\MM{X}}_{\beta}
\cdot\nabla\psi_k(\MM{X}_\beta), \\
& = & -\sum_lM_{kl}\langle \nabla\cdot([\tilde{\MM{u}}]\overline{D})\rangle_l.
\end{eqnarray*}
\qed

This means that is possible to take any reduced Lagrangian
$l(\MM{u},D)$ and, following the procedure described above, obtain a
finite-dimensional canonical Hamiltonian system with approximate
Lagrangian $\hat{L}\approx l$ for a set of particles, which may be
interpreted on the grid as a discrete form of the EP equations
resulting from $l$. If this system is integrated in time using a
symplectic method, then the spatially-discretised Hamiltonian will be
approximately preserved for long time-intervals. Suitable symplectic
methods are discussed in the next section.

\section{Symplectic time integration methods}
\label{time int}

In order to make a useful and practical Hamiltonian numerical method,
it is necessary to choose the symplectic time-integrator carefully.
In order to preserve the symplectic form, many integrators have to be
implicit; often the expense of solving the resulting algebraic
equations means that the symplectic integrator cannot compete with
non-conservative schemes. In this section we suggest a practical
first-order symplectic scheme for solving the semidiscretised system
obtained using the method of the previous section. In order to do
this, we restrict to a smaller set of problems where the Lagrangian
takes the form
$$
\mathcal{L} = \int \mathcal{K}(D,\MM{u})-\mathcal{V}(\MM{D})\diff^nx,
$$ where the nonlinear operator $\mathcal{K}$ is interpreted as the
kinetic energy and $\mathcal{V}$ is the potential energy; we shall
further require that $\mathcal{K}$ takes the form
$$
\mathcal{K} = \frac{1}{2}D\MM{u}\cdot\mathcal{B}\MM{u},
$$ where $\mathcal{B}$ is a linear operator, so that
$\MM{m}=\mathcal{B}^{-1}\MM{u}$. This set represents a
very large number of fluid systems so is not too much of a
restriction (in the examples we shall also consider the EPDiff
equation which does not have the density factor).
We shall proceed by discretising $\mathcal{L}$:
$$
\hat{L} = \sum_{kl}M_{kl}\left(
\frac{1}{2}D_l\MM{u}_l\cdot\sum_m(\hat{B}^{-1})_{lm}\MM{u}_m
-\hat{V}_l(\tilde{D})\right).
$$
For a given time step size $\Delta t$, write 
$$
\MM{X}_\beta^n \approx \MM{X}_\beta(n\Delta t),
$$
and similarly for $\overline{\MM{m}}$ \emph{etc}. To further compact 
the notation write
$$
[\tilde{f}]_\beta^n = \sum_k\tilde{f}_k\psi(\MM{X}_\beta^n),
\qquad
\langle\overline{f}\rangle_k^n = \sum_\beta\overline{f}_\beta
\psi(\MM{X}_\beta^n),
$$
and adopt similar notation for the gradient operators. We shall also
write $\MM{U}_k([\overline{\MM{m}}],\tilde{D})$ as the solution to
$$
\langle\overline{\MM{m}}\rangle_k = \nabla_{\MM{u}_k}\hat{L}
(\MM{U}(\langle\overline{\MM{m}}\rangle,\tilde{D})
,\tilde{D}),
$$
\emph{i.e.} the grid momentum-to-velocity map, so that
$$
\hat{K}_l = \frac{1}{2}\langle\overline{\MM{m}}\rangle_l\cdot
\MM{U}_l(\langle\overline{\MM{m}}\rangle,\tilde{D}).
$$

The (first-order) symplectic-Euler scheme then takes the form
\begin{eqnarray*}
\frac{\overline{\MM{m}}_\beta^{n+1}-\overline{\MM{m}}_\beta^n}
{\Delta t}
& = & - [(\nabla\MM{U}(\langle\overline{\MM{m}}^{n+1}\rangle^n
,\langle\overline{D}\rangle^n))^T]^n_\beta
\cdot\overline{\MM{m}}_\beta^{n+1} \\
& & \qquad -\overline{D}_\beta\left[\nabla\frac{1}{2}
\frac{\langle\overline{\MM{m}}^{n+1}\rangle^n}{\langle\overline{D}\rangle^n}
\cdot\MM{U}\left(
\langle\overline{\MM{m}}^{n+1}\rangle^n,\langle\overline{D}\rangle^n\right)
\right]^n_\beta \\
& & \qquad
-\overline{D}_\beta\left[\nabla\pp{V}{\tilde{D}}(\langle\overline{D}
\rangle^n)\right]^n_\beta, \\
\frac{\MM{X}_\beta^{n+1}-\MM{X}_\beta^n}
{\Delta t}
& = & - [\MM{U}(\langle\overline{\MM{m}}^{n+1}\rangle^n
,\tilde{D}^n)]^n_\beta,
\end{eqnarray*}
which is implicit in $\overline{\MM{m}}$ but not $\MM{X}$. The
algebraic system for $\overline{\MM{m}}$ can be efficiently solved
using the following iterative scheme:
\begin{eqnarray*}
\overline{\MM{m}}_\beta^{n+1,j+1}
& =& \overline{\MM{m}}_\beta^n - \Delta t[(\nabla\MM{U}(\langle\overline{\MM{m}}^{n+1,j}\rangle^n
,\langle\overline{D}\rangle^n))^T]^n_\beta
\cdot\overline{\MM{m}}_\beta^{n+1,j+1} \\
& & \qquad -\Delta t\overline{D}_\beta\left[\nabla\frac{1}{2}
\frac{\langle\overline{\MM{m}}^{n+1,j+1}\rangle^n}{\langle\overline{D}\rangle^n}
\cdot\MM{U}\left(
\langle\overline{\MM{m}}^{n+1,j}\rangle^n,\langle\overline{D}\rangle^n\right)
\right]^n_\beta \\
& & \qquad
-\Delta t\overline{D}_\beta\left[\nabla\pp{\hat{V}}{\tilde{D}}(\langle\overline{D}
\rangle^n)\right]^n_\beta,
\end{eqnarray*}
where $\overline{\MM{m}}_\beta^{n,j}$ is the $j$th element of the
sequence converging to $\overline{\MM{m}}_\beta^n$ (with the sequence
initialised with
$\overline{\MM{m}}_\beta^{n,0}=\overline{\MM{m}}_\beta^{n-1}$). This
results in inverting a sparse matrix (with number of non-zero elements
proportional to the number of particles) each iteration (and for the
case $\hat{V}=0$ the matrix becomes block diagonal as the different
particles decouple).

A second-order symplectic discretisation can be obtained by using the
second-order Lobatto IIIa-b partitioned Runge-Kutta method as
discussed in \cite{reic99:bac} with reference to adaptive symplectic
methods for molecular dynamics.

\section{Examples}
\label{examples}
\subsection{Example: the shallow-water-$\alpha$ equations}
\label{swalpha}
In this section we derive a Hamiltonian particle-mesh
semidiscretisation for the shallow-water-$\alpha$ (SW-$\alpha$)
equations. These equations are obtained by decomposing Lagrangian
particle trajectories into mean and fluctuating parts, averaging along
those trajectories, and then applying the ``frozen turbulence'' Taylor
assumption, before finally assuming isotropic, homogeneous turbulence
statistics \cite{hol99:_fluc_lag_eul}.

The model is derived from the reduced Lagrangian
$$
l = \int \frac{D}{2}\left(|\MM{u}|^2+\alpha^2|\nabla\MM{u}|^2
-gD\right)\diff{^2x},
$$
where $\alpha$ is the mean lengthscale for the fluctuations. The
variational derivatives are
$$
\dd{l}{\MM{u}} = (D+\nabla\cdot D\nabla)\MM{u},
\qquad \dd{l}{D} = \frac{1}{2}\left(|\MM{u}|^2+\alpha^2|\nabla\MM{u}|^2\right)-gD,
$$
and substituting these derivatives into the EP equation
$$
\partial_t\MM{m}+(\nabla\MM{u})^T\cdot\dd{l}{\MM{u}}
+\nabla\cdot\left(\MM{u}\dd{l}{\MM{u}}\right) = D\nabla\dd{l}{D},
$$
gives the equations
\begin{eqnarray}
\partial_t\MM{m} + \MM{u}\cdot\nabla\MM{m}
+(\nabla\MM{u})^T\cdot\MM{m} + \MM{m}\nabla\cdot\MM{u}
& = &
D\nabla\Big(-g D + \nonumber \\
& & \qquad \frac{1}{2}\left(|\MM{u}|^2+\alpha^2|\nabla\MM{u}|^2\right)
\Big), \nonumber \\
D_t+\nabla\cdot(D\MM{u}) & = & 0, \nonumber \\
\MM{m} &=& (D-\alpha^2\nabla\cdot D\nabla)\MM{u}. \label{elliptic}
\end{eqnarray}
In these equations the momentum $\MM{m}$ is advected by the velocity
$\MM{u}$ obtained by inverting equation (\ref{elliptic}). The
equations have a more familiar appearance in the form
\begin{multline*}
(\partial_t+\MM{u}\cdot\nabla)(1-\alpha^2D^{-1}\nabla\cdot D)^{-1}
\MM{u}
+(\nabla\MM{u})^T\cdot(1-\alpha^2D^{-1}\nabla\cdot D)^{-1}\MM{u} \\
\shoveright{ = \nabla\left(-g D + \frac{1}{2}\left(|\MM{u}|^2+\alpha^2|\nabla\MM{u}|^2\right)
\right),} \\
 \qquad D_t+\nabla\cdot(\MM{u}D)=0,
\end{multline*}
which become the shallow-water equations
$$
(\partial_t+\MM{u}\cdot\nabla)\MM{u} = -g\nabla D, \qquad
D_t+\nabla\cdot(\MM{u}D)=0,
$$
when $\alpha^2$ is set to zero, making use of the identity
$$
\nabla\frac{1}{2}|\MM{u}|^2 = (\nabla\MM{u})\cdot\MM{u}.
$$

Following the procedure set out in section \ref{approach},
we discretise the Lagrangian using a piecewise-linear Galerkin finite element
discretisation
(see \cite{red93:int_fin_elem_met} for example). The
Galerkin expansion for a function $f$ given on grid points is
$$
f(\MM{x}) = \sum_kN_k(\MM{x})f_k,
$$
where $N_k$ is a continuous basis function which satisfies
$$
N_i(\MM{x}_j) = \delta_{ij}.
$$
Under this scheme, the Lagrangian becomes
$$
\hat{L} = 
\sum_{kl}\left(\frac{1}{2}\tilde{\MM{u}}_k\cdot B_{kl}\tilde{\MM{u}}_k
-\frac{1}{2}g\tilde{D}_kM_{kl}\tilde{D}_l
\right),
$$
where the Helmholtz matrix $B$ takes the form
$$
{B}_{ij}=\int_{\Omega}\sum_k\tilde{D}_kN_k(\MM{x})
\left(N_i(\MM{x})N_j(\MM{x})+\alpha^2\nabla N_i(\MM{x})\cdot
\nabla N_j(\MM{x})\right)\diff^2x,
$$
and the mass matrix $M$ takes the form
\begin{equation}
\label{mass}
M = \int_{\Omega}N_i(\MM{x})N_j(\MM{x})\diff^2x.
\end{equation} 
Both $\hat{B}$ and $M$ are sparse, well-conditioned, symmetric
matrices so they are very quick to invert using a preconditioned
conjugate gradients algorithm.

Using the Legendre transform, the momentum $\tilde{m}_k$ on the grid
is 
$$
\tilde{\MM{m}}_k = \pp{\hat{L}}{\tilde{\MM{u}}_k}=
\sum_lB_{kl}\tilde{\MM{u}}_k,
$$
and the discrete Hamiltonian takes the form
$$
\hat{H} = \sum_{kl}\left(\frac{1}{2}\tilde{\MM{m}}_k
  \cdot(M_{kl}\sum_mB^{-1})_{lm}\tilde{\MM{m}}_m
  +\frac{1}{2}g\tilde{D}_kM_{kl}\tilde{D}_l
\right).
$$

The canonical equations derived from the Hamiltonian $\hat{H}$ are
then
\begin{eqnarray*}
\dot{\MM{X}}_{\beta} &=& [\tilde{u}]_\beta=
[B^{-1}\langle\overline{\MM{m}}\rangle]_\beta, \\
\frac{\dot{\overline{\MM{m}}}_\beta}{D_\beta} 
&=& -[(\nabla\tilde{\MM{u}})^T]_\beta\cdot
\frac{\overline{\MM{m}}_\beta}{D_\beta} -
g[\nabla\tilde{D}]_\beta-\left[\nabla\frac{1}{2}
\frac{\langle\overline{\MM{m}}\rangle}{\langle\overline{D}\rangle}
\cdot B^{-1}
\langle\overline{\MM{m}}\rangle
\right]_\beta,
\end{eqnarray*}
and the momentum equation on the grid take the form
$$
\frac{d}{dt}\langle\overline{\MM{m}}\rangle_k
+\left\langle\nabla\cdot\left(\left[\tilde{\MM{u}}\right]
\overline{\MM{m}}\right)\right\rangle_k
+\left\langle
\left[(\nabla\tilde{\MM{u}})^T\right]
\cdot
\overline{\MM{m}}
\right\rangle_k
=-\left\langle \overline{D}
\left[\nabla\left(gD+\frac{1}{2}
\frac{\langle\overline{\MM{m}}\rangle}{\langle\overline{D}\rangle}\cdot
\tilde{\MM{u}}\right)
\right]\right\rangle_k,
$$

We can also make a semidiscretisation for the Euler-$\alpha$ equations
by removing the potential energy term and setting a constraint
$$
D(\MM{x}_k)=1, \qquad k=1,\ldots,m,
$$ \emph{via} a set of Lagrange multipliers which give the pressure at
each grid-point. It is possible to construct symplectic integrators
for constrained Hamiltonian systems of this type
\cite{jay96:sym_run_kut_ham,leim94:sym_ham,reic96:sym_ham} and these
methods have already been applied to enforcing an incompressibility
constraint in \cite{cot03:hamiltonian_particle_mesh_method_two}.

\subsection{Example: the 2D EP-Diff equations}
\label{epdiff}
The EP-Diff equations \cite{hol05:mom_epd} are the generalisation to
$n$-dimensional space of the Camassa-Holm equation \cite{cam93}, and
represent the geodesic equations of motion for the $H^1$ norm. They
represent the prototype for the family of regularised equations such
as Euler-$\alpha$ \cite{hol99:_fluc_lag_eul,mar00:eul}, although they
also have applications in computational anatomy and image processing
\cite{mil02:eul_lag,mum98:ques_mat_en_trait_du_sig}.

The 2D EP-Diff equation is obtained from the reduced Lagrangian:
$$
l = \int \frac{|\MM{u}|^2}{2}+\alpha^2\frac{|\nabla\MM{u}|^2}{2}\diff{^2x},
$$ where $\alpha$ is a constant lengthscale parameter, and so the
equation of motion is
$$
\partial_t\MM{m}+\MM{u}\cdot\nabla\MM{m}
+(\nabla\MM{u})^T\cdot\MM{m}+\MM{m}\nabla\cdot\MM{u}=0,
\qquad \MM{m} = (1-\alpha^2\Delta)\MM{u}.
$$

Following the programme set out in section \ref{approach} with a
Galerkin finite element approximation, the discrete
Lagrangian takes the form
$$
\hat{L} = \sum_{kl}\frac{1}{2}\tilde{\MM{u}_k}
\cdot (A_{kl}\tilde{\MM{u}}_l).
$$
where the matrix $A$ is the discrete inverse
modified Helmholtz operator 
$$
A_{kl} = \int_{\Omega}N_k(\MM{x})N_l(\MM{x})+\alpha^2
\nabla N_k(\MM{x})\cdot\nabla N_l(\MM{x})\diff^2x.
$$

The discrete Hamiltonian in
Eulerian coordinates on the grid is then:
$$
H = \sum_{kl}\frac{1}{2}\tilde{\MM{m}}_k
\cdot\sum_m(A^{-1})_{km}M_{ml}\tilde{\MM{m}}_l,
$$
where $M_{kl}$ is the finite element mass matrix given by
equation (\ref{mass}).

The canonical equations for this Hamiltonian are then
\begin{eqnarray*}
\dot{\MM{X}}_\beta & = & [\MM{U}(\langle\overline{\MM{m}}\rangle)
]_{\beta} = 
[A^{-1}\langle\overline{\MM{m}}\rangle]_\beta, \\
\frac{\dot{\overline{\MM{m}}}_\beta}{\overline{D}_\beta}
& = & -[(\nabla\MM{U}(\langle\overline{\MM{m}}\rangle))^T]_\beta\cdot\MM{m}_\beta,
\end{eqnarray*}
and the equation for $\MM{m}$ is the EP-Diff equations represented
in our discrete calculus.

The symplectic Euler method with these equations is 
\begin{eqnarray*}
\frac{\overline{\MM{m}}_\beta^{n+1}-\overline{\MM{m}}_\beta^n}
{\Delta t}
& = & - [(\nabla\MM{U}(\langle\overline{\MM{m}}^{n+1}\rangle^n))^T
]^n_\beta
\cdot\overline{\MM{m}}_\beta^{n+1} \\
\frac{\MM{X}_\beta^{n+1}-\MM{X}_\beta^n}{\Delta t}& = & [\MM{U}(\langle\overline{\MM{m}}^{n+1}\rangle^n)]^n_\beta,
\end{eqnarray*}
and the iterative scheme to solve for $\overline{\MM{m}}_\beta^{n+1}$ is
\begin{equation}
\label{epdiff iterative}
\overline{\MM{m}}_\beta^{n+1,j+1}=\overline{\MM{m}}_\beta^n - \Delta t[(\nabla\MM{U}(\langle\overline{\MM{m}}^{n+1,j}\rangle^n))^T]^n_\beta
\cdot\overline{\MM{m}}_\beta^{n+1,j+1}.
\end{equation}
Using this scheme only requires the inversion of one $3\times 3$ matrix
for each particle per iteration which leads to a very efficient method.

\section{Numerical results}
\label{numerics}
To demonstrate the particle-mesh approach, we give results using the
method for the 2D EP-Diff equation described section \ref{epdiff}
(numerical results and discussion of emergent behaviour for this
system are given in \cite{hol03:wav_pdes}). The equations are solved
in a $2\pi\times 2\pi$ periodic square, with $\alpha=0.3133$. The
discrete Lagrangian was obtained using piecewise-linear quadrilateral
finite elements on a regular $128\times 128$ grid and the particles
positions were initialised 16 particles in each grid cell. The
equations were integrated using Matlab using the HPM C-mex files (see
{\ttfamily http://www.cwi.nl/projects/gi/HPM/}) to perform the
grid-particle operations.  The finite element matrices $B$ and $M$
were stored as sparse matrices and inverted using conjugate gradients
with incomplete Cholesky preconditioning, reaching a residual of less
that $1\times 10^{-9}$ in less than 10 iterations.  The iterative
scheme (\ref{epdiff iterative}) converged with a total error of less
that $1\times 10^{-9}$ in less than 3 iterations for each timestep
(with $dt=0.0204$) during the experiment. 

\begin{figure}[h]
\begin{center}
\scalebox{0.75}{\includegraphics{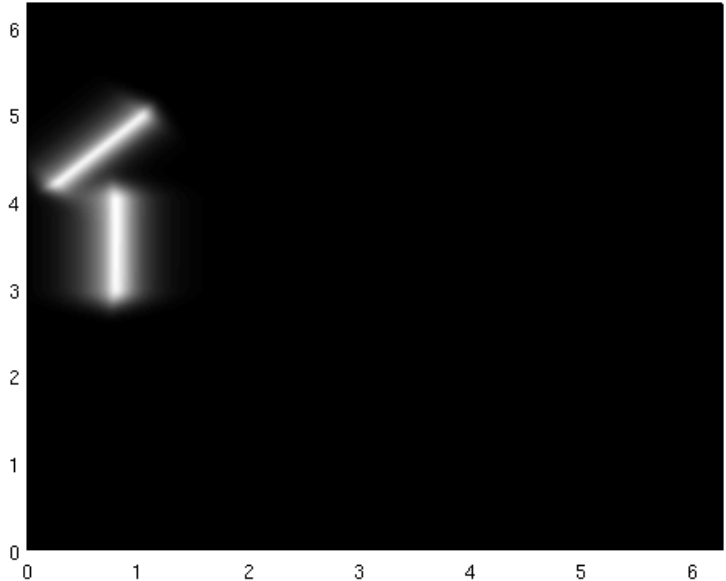}}
\end{center}
\caption{\label{001}Plot showing magnitude of velocity at $t=0$. The
  velocity is distributed along two lines which are set to collide.}
\end{figure}

\begin{figure}[h]
\begin{center}
\scalebox{0.75}{\includegraphics{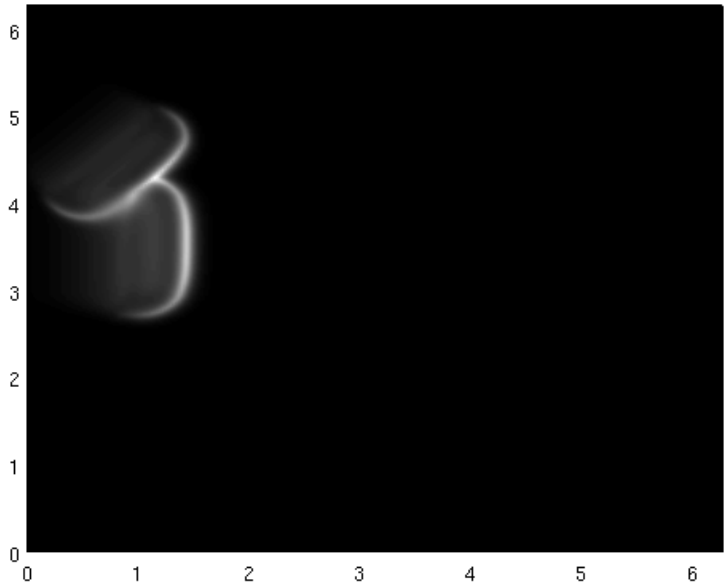}}
\end{center}
\caption{\label{012}Plot showing magnitude of velocity at
  $t=0.2445$. The two lines have stretched out and expanded before overlapping.} 
\end{figure}

\begin{figure}[h]
\begin{center}
\scalebox{0.75}{\includegraphics{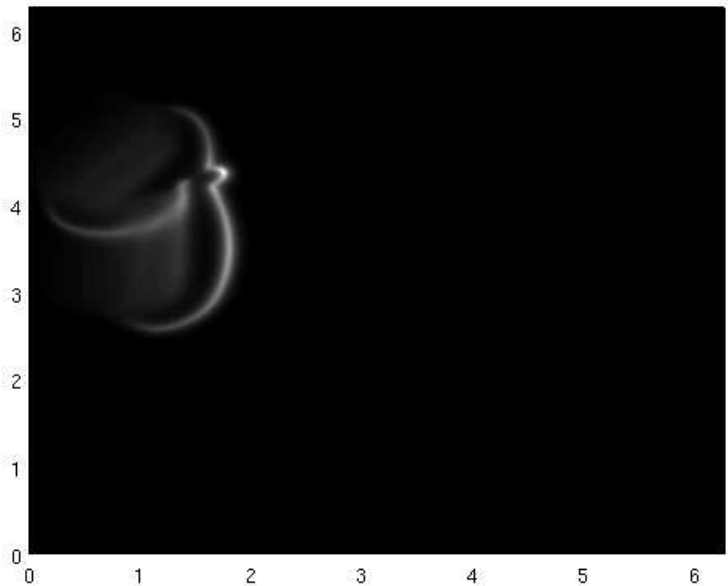}}
\end{center}
\caption{\label{020}Plot showing magnitude of velocity at
  $t=0.4074$. The lines have reconnected after crossing.}
\end{figure}

\begin{figure}[h]
\begin{center}
\scalebox{0.75}{\includegraphics{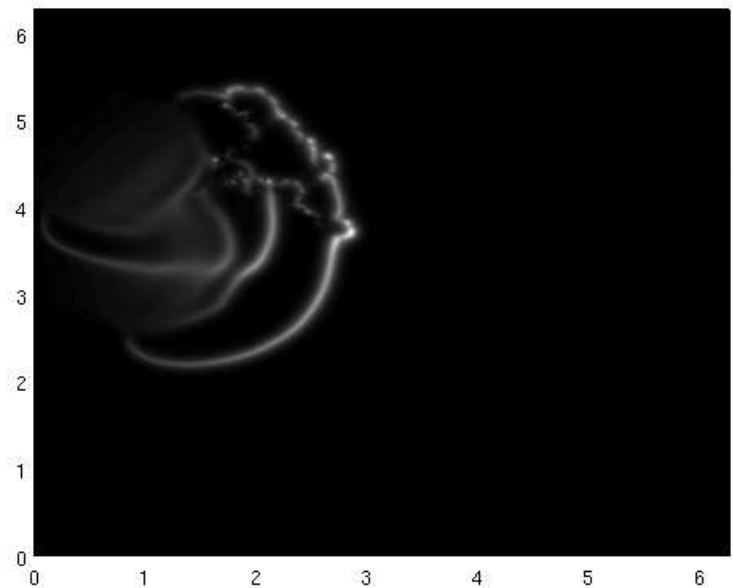}}
\end{center}
\caption{\label{045}Plot of showing magnitude of velocity at
  $t=0.9167$. The lines have separated into a series of thinner lines
  which are the emergent structures for these equations.}
\end{figure}

\begin{figure}[h]
\begin{center}
\scalebox{0.75}{\includegraphics{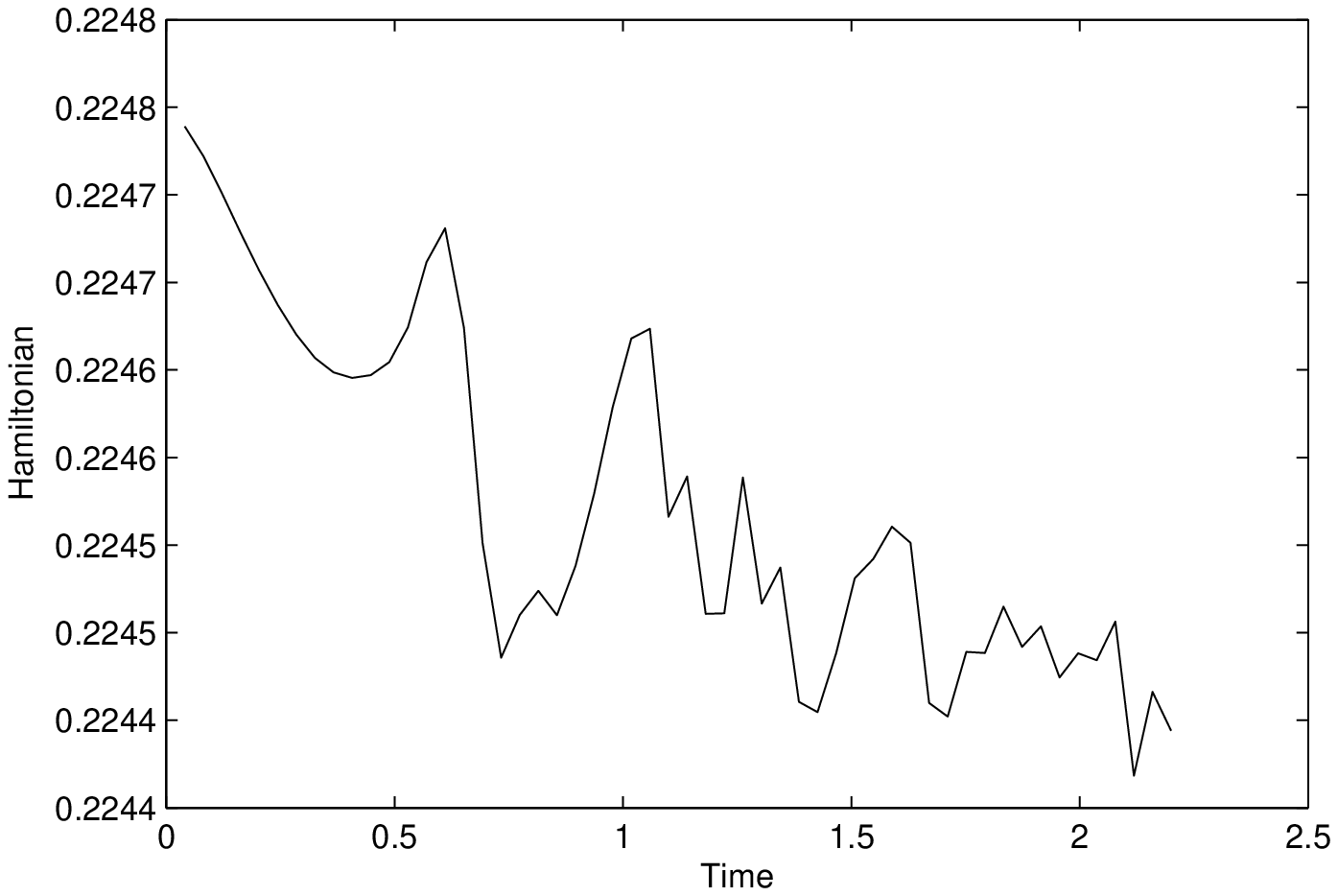}}
\end{center}
\caption{\label{Ham}Plot of the Hamiltonian against time during the numerical experiment. The time stepping method is symplectic and so the Hamiltonian stays within $\mathcal{O}(\Delta t)$ of the value at $t=0$ for very long time-intervals.}
\end{figure}

The initial data had the velocity zero everywhere except for two thin
lines which are set to collide. The plots of $x$-velocity in figures
\ref{001}-\ref{045} show the lines reconnecting after collision and
the formation of thinner ``peak-shaped'' lines, as reported in
\cite{hol03:wav_pdes}. A plot of the discrete Hamiltonian during the
numerical experiment is shown in figure \ref{Ham}. This plot
illustrates that, as the timestepping method used is symplectic, the
discrete Hamiltonian is conserved within $\mathcal{O}(\Delta t)$ of
the initial value at time $t=0$ for very long time-intervals. This is
the real benefit of discretising the variational structure of the
equations. Higher-order time-integrators will produce smaller errors
in the Hamiltonian as $\Delta t\to0$.

\section{Discussion and outlook}
\label{discussion}
We introduced a general framework for deriving canonical Hamiltonian
semi-discretised equations from the reduced Lagrangian in Eulerian
coordinates, based on the particle-mesh approach. This framework has
several advantages:
\begin{itemize}
\item The use of Lagrangian particles in the discretisation means that
  the equations are always canonical. This means that the discrete
  structure operator is guaranteed to satisfy the Jacobi identity and
  so a symplectic time integrator can be used, leading to long-time
  preservation of the Hamiltonian. The canonical structure also makes
  for easier discussion of adiabatic invariants in the discrete system
  \cite{cot04:mod}.
\item The Hamiltonian is given in terms of Eulerian coordinates. This
  means it is easy to add on extra physics, make approximations,
  \emph{etc}. Many models, such as those discussed in section
  \ref{examples}, involve differential operators given in Eulerian
  coordinates which become rather complicated and ``entangled'' in the
  Lagrangian variables.  This approach avoids this problem.
\item This approach allows the momentum densities to become
  $\delta$-functions in the limit, whilst the velocities on the grid
  can be represented using finite elements which means that weak
  solutions can be discussed, for example in the EP-Diff equations we
  seek solutions in $H^1$ (as illustrated in the numerical example).
\item The symplectic Euler method can be efficiently applied to
  integrate the spatially-discrete equations in space. It should also
  be possible to obtain higher-order symplectic methods for greater
  accuracy.
\item The numerical experiment demonstrates that, at least for the
  EP-Diff equations, this numerical method is computationally feasible
  and produces well-behaved numerical solutions.
\end{itemize}
Future work may take several different directions, including:
\begin{itemize}
\item Development of higher-order symplectic integrators for these models
  in order to produce practical methods.
\item Application to other fluid models: the two-layer Green-Naghdi
  equations \cite{choi99:ful}, the quasi-geostrophic equations
  \cite{hol98}, rotating shallow-water-$\alpha$, Euler-Boussinesq-$\alpha$
  in 3D \cite{hol99:_fluc_lag_eul} \emph{etc}.
\item Comparison of methods with existing schemes for simple models
  such as the EP-Diff equations in 1, 2 and 3 dimensions.
\item Investigation of the effects of smoothing in the
  shallow-water-$\alpha$ equations. Does the smoothing in the velocity 
  field keep the solution well-behaved over long time intervals?
\item Use of normal-form theory, following the work in
 \cite{cot04:mod}, to investigate geostrophic balance in the numerical
 scheme for the rotating shallow-water-$\alpha$ equations, and other
 rotating models.
\item Comparison of the solutions of shallow-water-$\alpha$ obtained
  using the method given in this paper with solutions of RLDSW
  obtained using the HPM method.
\item Inclusion of boundary conditions: the approach, as described
  here, only works on a closed manifold (such as the 2-torus
  \emph{i.e.}  a square with periodic boundary conditions).
\item Obtain an extension to semi-direct product systems with other
  types of advected quantities (other than densities) such as scalars,
  vectors, and higher-order tensors.
\end{itemize}

\noindent {\bfseries Acknowledgements:} Thanks to Sebastian Reich,
Darryl Holm and Joel Fine for their useful discussions during the
writing of this paper, to Jason Frank for the use of his HPM C-mex
code, to Matthew Piggott for advice on programming the finite element
method, to the Daniels family of Ohio for their generous hospitality
during the early stages of this work, and to NERC for funding under
grant number NER/T/S/2002/00459.

\nocite{*}

\bibliography{generalPM}

\end{document}